\pdfoutput=1
\documentclass[12pt]{article}

\usepackage{amssymb}

\usepackage{alltt,times,pifont,amssymb,url,graphicx} % JRH \{jrh}
\usepackage{latexsym} % RDA

\def\Hide#1{\relax}

\DeclareSymbolFont{AMSb}{U}{msb}{m}{n}
\DeclareSymbolFontAlphabet{\mathbb}{AMSb}

\def\Func#1{{\mathsf{#1}}}

\def\fin{\mathbb{F}}

\def\VecSort{{\cal V}}
\def\ScaSort{{\cal R}}

\def\NS#1{\Func{NS}^{#1}}

\def\MS{\Func{MS}}
\def\NSA#1{\NS{#1}_{{+}}}

\def\MSA{\Func{MS}_{{+}}}

\def\lL{{\cal L}}

\def\LM{\lL_{M}}
\def\LN{\lL_{N}}
\def\LPA{\lL_{PA}}
\def\LNA{\lL^{+}_{N}}

\def\sS{\mathbb{S}}
\def\sJ{\mathbb{J}}

\def\sV{\mathbb{V}}

\def\pADS#1{\Func{ADS}(#1)}

\def\pEP#1{\Func{EP}(#1)}

\def\pM#1{\Func{M}(#1)}
\def\pMGI#1{\Func{M}_{{>}1}(#1)}

\def\pOK{\Func{OK}}

\def\pS{\Func{S}}
\def\pHPV#1{\Func{HPV}(#1)}
\def\pHPL#1{\Func{HPL}(#1)}
\def\pSEP#1{\Func{SEP}(#1)}

\def\sS{\Func{S}}

\def\IsDef{\mathrel{{:}{=}}}

\def\Imp{\Rightarrow}

\def\And{\land}

\def\ST{\mathrel{|}}

\def\To{\rightarrow}

\def\Vb{\mathbf{b}}

\def\Ve{\mathbf{e}}

\def\Vp{\mathbf{p}}
\def\Vq{\mathbf{q}}
\def\Vr{\mathbf{r}}

\def\Vu{\mathbf{u}}
\def\Vv{\mathbf{v}}
\def\Vw{\mathbf{w}}
\def\Vx{\mathbf{x}}
\def\Vy{\mathbf{y}}

\def\VO{\mathbf{0}}

\def\Norm#1{{||}#1{||}}

\newtheorem{Theorem}{Theorem}
\newtheorem{Lemma}[Theorem]{Lemma}
\newtheorem{Corollary}[Theorem]{Corollary}
\def\Proof{\par \noindent{\bf Proof: }}
\def\Example{\par \noindent{\bf Example: }}
\def\Done{\hfill\rule{0.5em}{0.5em}}

\newcommand{\nat}{\mbox{$\protect\mathbb N$}}
\newcommand{\natgtI}{\mbox{${\protect\mathbb N}_{{>}1}$}}

\newcommand{\rat}{\mbox{$\protect\mathbb Q$}}
\newcommand{\ratgtO}{\mbox{${\protect\mathbb Q}_{{>}0}$}}
\newcommand{\real}{\mbox{$\protect\mathbb R$}}

\newcommand{\spot}{{\cdot}}

\newcommand{\all}[1]{\forall #1 \spot\:}
\newcommand{\ex}[1]{\exists #1 \spot\:}
\newcommand{\exu}[1]{\exists! #1 \spot\:}

\newcommand{\Ands}{\bigwedge}

\newcommand{\BEQ}{\mbox{\raise4pt\hbox{$\ulcorner$}}}

\newcommand{\EEQ}{\mbox{\raise4pt\hbox{$\urcorner$}}}

\newcommand{\BA}{\begin{array}[t]{l}}
\newcommand{\EA}{\end{array}}

\makeatletter
\def\imod#1{\allowbreak\mkern10mu({\operator@font mod}\,\,#1)}
\makeatother

\newcommand\HOLSpacing{12pt}

\newlength{\hsbw}
\title{The decision problem for normed spaces over any class of ordered fields}
\author{R.D. Arthan\thanks{Inspired by joint
work with Robert M. Solovay and John Harrison.}\\
Lemma 1 Ltd, 2nd Floor, 31A Chain Street, Reading, RG1 2HX\\
E-mail: {\tt rda@lemma-one.com}}

\date{April 17, 2011}

\begin{document}
\maketitle

\begin{abstract}
It is known that the theory of any class of normed spaces over $\real$ that
includes all spaces of a given dimension $d \ge 2$ is undecidable, and indeed,
admits a relative interpretation of second-order arithmetic.  The notion of a
normed space makes sense over any ordered field of scalars, but such a strong
undecidability result cannot hold in the more general case. Nonetheless, we
find that the theory of any class of normed spaces in the more general sense
that includes all spaces of a given dimension $d \ge 2$ over some ordered field
admits a relative interpretation of Robinson's theory $Q$ and hence is
undecidable.
\end{abstract}

Let $\LN$ be the natural two-sorted language for a normed space over an
arbitrary ordered field of scalars:  $\LN$ has sorts $\cal K$ for the
scalars and $\VecSort$ for the vectors together with  the usual
symbols of the appropriate sorts for a vector space over an ordered
field equipped with a norm
(see \cite{Solovay-et-al09} for more details; we have adopted $\cal K$
instead of $\ScaSort$ for the scalar sort here).
If $K$ is an ordered field, a {\em normed space over $K$} is a structure
for $\LN$ in which $\cal K$ and the field symbols are interpreted
in $K$ and which satisfies the usual first-order axioms for a normed space.
If $\cal C$ is any non-empty class of ordered
fields, let $\NS{\relax}({\cal C})$ be the set of all sentences in $\LN$
that are valid in all normed spaces whose field of scalars belongs to $\cal C$;
let 
$\NS{n}({\cal C})$ for $n \in \nat$, 
$\NS{\fin}({\cal C})$
and $\NS{\infty}({\cal C})$
denote the extensions of the theory $\NS{\relax}({\cal C})$ comprising
the sentences
valid in all normed spaces whose field of scalars belongs to $\cal C$
and that are,  respectively,
$n$-dimensional, finite-dimensional and infinite-dimensional;
let 
$\NSA{n}({\cal C})$,
$\NSA{\fin}({\cal C})$
and $\NSA{\infty}({\cal C})$ be the corresponding theories in the purely
additive sublanguage $\LNA$, i.e., the sublanguage in which
scalar-scalar and scalar-vector multiplication are disallowed (although
we may still use multiplication by rational constants as a shorthand).

\cite{Solovay-et-al09} deals with the case ${\cal C} = \{\real\}$ and shows
that with the exception of the 1-dimensional case (which
reduces trivially to the first-order theory of $\real$),
all of the theories mentioned in the previous paragraph are
undecidable with this choice of $\cal C$ in the strong sense that they admit a relative
interpretation of second-order arithmetic.
In general, $\cal C$ may be definable by a recursive set
of axioms (the class of real closed fields is an example).  In this case,
$\NS{\relax}({\cal C})$ is recursively axiomatizable,
implying that we cannot interpret second-order arithmetic in it. Nonetheless, we shall see
that for any non-empty $\cal C$, even the purely
additive theory $\NSA{\relax}({\cal C})$ is undecidable, as are
all of $\NSA{n}({\cal C})$ for $1 < n \in \nat$, 
$\NSA{\fin}({\cal C})$ and $\NSA{\infty}({\cal C})$.

We will use the classical method of proving that a theory $T$ is undecidable
by giving a relative interpretation in $T$ of Raphael M. Robinson's finitely
axiomatizable and essentially undecidable theory $Q$. Recall, e.g., from
\cite{tarski-et-al53}, that $Q$ is the theory in the language $\LPA$
of Peano arithmetic comprising the deductive closure of the following
axioms:
$$
\begin{array}{lc}
\mbox{Q1:} & \all{x\;y} \sS(x) = \sS(y) \Imp x = y \\
\mbox{Q2:} & \all{x} 0 \not= \sS(x) \\
\mbox{Q3:} & \all{x} x \not= 0 \Imp \ex{y} x = \sS(y) \\
\mbox{Q4:} & \all{x} x + 0 = x \\
\mbox{Q5:} & \all{x\;y} x + \sS(y) = \sS(x + y) \\
\mbox{Q6:} & \all{x} x \times 0 = 0 \\
\mbox{Q7:} & \all{x\;y} x \times \sS(y) = x \times y + x
\end{array}
$$

In $\LNA$, the intended interpretation of $\cal K$ is as an ordered group $K$
with a distinguished positive element $1$.  If $K$ is such a group, let us
write $\nat_K$ for the semiring of natural numbers considered as a subset of
$K$ by identifying $n$ with $\sum_{i=1}^n1$.  We then have the following
theorem:

\begin{Theorem}\label{thm:gen-linear-undec}
Let $L$ be a (many-sorted) first-order language including a sort $\cal K$,
together with a function symbol ${+} : {\cal K} \times {\cal K} \rightarrow
{\cal K}$, a binary predicate symbol ${<}$ on the sort $\cal K$ and a
constant symbol $1 : \cal K$ whose intended interpretations are as some
ordered abelian group with $1$ as a distinguished positive element.  Let $\cal
C$ be some class of structures for $L$, in which $\cal K$ and these symbols
have their intended interpretations and let $\cal T$ be the theory of $\cal C$,
i.e., the set of all sentences of $L$ valid in every member of $\cal C$.  Let $\mu(x,
y, z)$ be a formula of $L$ with the indicated free variables all of sort $\cal
K$.  Let $\cal M$ be a structure in the class $\cal C$, in which $\cal K$ is
interpreted as some ordered group $K$.  and assume that in $\cal M$, $\mu(x, y,
z)$ defines the graph of multiplication in $\nat_K$. Then $\cal T$ is
undecidable.
\end{Theorem}
{\Proof}
Define $\nu(x) \IsDef \mu(x, 0, 0)$, so that in $\cal M$, $\nu(x)$ holds iff $x \in
\nat_K$. 
Define a translation $\phi \mapsto \phi^*$ from sentences of
$\LPA$ to sentences of $L$, where $\phi^*$ is obtained
from $\phi$ by the following sequence of transformations:
{\em(i)} label all constants and variables with the sort $\cal K$;
{\em(ii)} unnest ocurrences of $\times$ so that
$\times$ only occurs in formulas of the form $z = x\times y$
where $x$, $y$ and $z$ do not involve $\times$;
{\em(iii)} relativise with
respect to $\nu$, i.e., replace all subformulas of the form $\all{x}\psi$
(resp. $\ex{x}\psi$) by $\all{x}\nu(x) \Imp \psi$ (resp. $\ex{x}\nu(x) \And
\psi$);
{\em(iv)} replace all subterms of the form  $\pS(x)$ by $x + 1$;
and {\em(v)} replace subformulas of the form 
$z = x \times y$ by $\mu(x, y, z)$.

Define a sentence $\pOK$ of $L$ as follows:
\begin{eqnarray*}
\pOK &\IsDef&
 \BA
   \nu(0) \And {} \\
   (\all{x} \nu(x) \Imp x \ge 0 \And \nu(x + 1)) \And {} \\
   (\all{x} \nu(x) \And x > 0  \Imp \nu(x - 1)) \And {} \\
   (\all{x\;y}\nu(x) \And \nu(y) \Imp \exu{z} \mu(x,y,z)) \And {} \\
   (\all{x\;y\;z} \mu(x, y, z) \Imp \nu(x) \And \nu(y) \And \nu(z) ) \And {} \\
   (\all{x\;y\;w\;z} \mu(x, y + 1, w) \And \mu(x, y, z) \Imp w = z + x)
  \EA
\end{eqnarray*}
and write $\pOK_i$ for the $i$-th conjunct in $\pOK$.

It is easy to verify that the translations $\mbox{Q1}^* \ldots \mbox{Q7}^*$ of
the axioms of $Q$ all hold in any normed space over any ordered field in which
$\pOK$ holds.  For example, $\mbox{Q6}^*$ is equivalent to the tautology
$\all{x}\nu(x) \Imp \nu(x)$ and $\mbox{Q3}^*$ holds because if $x \not= 0$ and
$\nu(x)$ holds, then, by $\pOK_2$, $x \ge 0$, whence $x > 0$, so that
$\nu(x-1)$ holds by $\pOK_3$, so that $(\ex{y}x = \sS(y))^*$ holds with $x-1$
as witness.  We have constructed a relative interpretation of the essentially
undecidable and finitely axiomatizable theory $Q$ in the theory, $T_1$ say,
obtained by adding the finite set of axioms $\{\pOK\}$ to $T$.  Clearly $\cal
M$ is a model of $T_1$, so $T_1$ is consistent.  It follows from Theorems I.8
and I.10 of \cite{tarski-et-al53} that $T$ is undecidable.
\Done

\Example The following is based on an idea of John Harrison.
If $K$ is any ordered field, define a {\em metric space over $K$} to be a set
$X$ equipped with a function $d : X \times X \To K$ satisfying the usual axioms
for a metric space.  Now assume that $K$ is a euclidean field, i.e., an ordered
field such that every positive element has a square root, so that the vector
space $K^n$ admits the
euclidean norm $\Norm{(x_1, \ldots, x_n)}_e = \sqrt{x_1^2+ \ldots + x_n^2}$
and becomes a metric space over $K$ under $d_e(\Vx, \Vy) = \Norm{\Vy - \Vx}_e$.
Let $\sV = X \cup Y \subseteq K^2$ where $X =
\{(x, y) | x \in \nat_K \And y = 0\}$ and $Y = \{(x, y) | x \in \nat_K \And y =
x^2\}$.  In $K^n$ just as in $\real^n$, a
point $\Vq$ lies on the line segment  $[\Vp, \Vr]$ iff $d_e(\Vp, \Vr) =
d_e(\Vp, \Vq) + d_e(\Vq, \Vr)$.  Using this fact, it is not difficult to give a
first-order formula $\phi(x, y)$ in the natural two-sorted language $\LM$ for
metric spaces with the indicated free variables of sort $\cal K$ that holds
in $\sV$ iff $x \in \nat_K$ and $y = x^2$ (Design $\phi(x, y)$ to assert that
$x = d(\VO, \Vx)$ and $y = d(\Vx, \Vy)$ where
$\Vy \in Y$ and $\Vx$ is the point of $X$ nearest to $\Vy$.  Cf. the proof of
Theorem~5 in \cite{Solovay-et-al09}.).  Now if we put $\mu(x, y, z) \equiv
\ex{a\;b\;c} \phi(x, a) \And \phi(y, b) \And \phi(x + y, c) \And
a + b + 2z = c$, then $\mu(x, y, z)$ defines the graph of multiplication in
$\nat_K$.  Applying the theorem, we obtain the undecidability of the theory
$\MS({\cal C})$ of metric spaces over any class $\cal C$ of ordered fields that
includes at least one euclidean field. (In fact, $\phi$ can be defined without
using multiplication, so the additive theories $\MSA({\cal C})$ are also
undecidable.)
\Done

We will give a construction inspired by the proof of theorem~41 in
\cite{Solovay-et-al09}, where we found normed spaces over $\real$ in which
there are definable consecutive pairs of line segments inscribed in the unit
circle whose lengths are in the ratio $1:m$ for $m$ in the set $\natgtI$ of
natural numbers greater than 1.
Now, working over an arbitrary ordered field
$K$, we will construct normed spaces $\sJ^d$ in which there are
definable consecutive quadruples of line segments inscribed in the unit circle
$S$ whose lengths are in the ratio $1:m:mn:n$ for $m, n \in \natgtI$. Thus for
positive $r \in K$ if one of the corresponding quadruples $(x_1, x_2, x_3,
x_4)$ in the circle $rS$ has $x_1 = 1$ then $x_2$ and $x_4$ are in $\nat_K$ 
and $x_3 = x_2x_4$; moreover, for any $x_2, x_4 \in \nat_K$, such a
quadruple exists for some $r > 0$. This will allow us to apply
theorem~\ref{thm:gen-linear-undec} to conclude the undecidability of any class
of normed spaces that includes $\sJ^d$.

The proof of theorem~41 in \cite{Solovay-et-al09}  was based on the convergence of the
power series for the exponential function. We need a replacement for this
series that has consecutive quadruples of terms in the ratios
$1:m:mn:n$ for all $m , n\in \natgtI$.
The following lemma gives us this.

\begin{Lemma}\label{lma:a-etc}
There are $p, q, m_i, n_i \in \natgtI$ and 
$a, a_k \in \ratgtO$, $i = 0, 1, \ldots$, $k = 1, 2, \ldots$
satisfying the following conditions:
\begin{description}
\item[{\em(i)}] the function $(m, n) \mapsto p^mq^n$ is an injection of $\natgtI \times \natgtI$ into $\natgtI$;
\item[{\em(ii)}] the pairs $(m_i, n_i)$, $i = 0, 1, \ldots$ enumerate $\natgtI \times \natgtI$;
\item[{\em(iii)}] $ p^{m_0}q^{n_0} < p^{m_1}q^{n_1} < p^{m_2}q^{n_2} < \ldots $;
\item[{\em(iv)}] for $k = 1, 2, \ldots$,
\begin{eqnarray*}
a_k &\IsDef&
\left\{
    \begin{array}{rl}
        p^{-m_i}q^{-n_i} &  \mbox{if $k = 4i + 1$} \\
        m_i p^{-m_i}q^{-n_i} & \mbox{if $k = 4i + 2$} \\
        m_i n_i p^{-m_i}q^{-n_i} & \mbox{if $k = 4i + 3$} \\
        n_i p^{-m_i}q^{-n_i} & \mbox{if $k = 4i + 4$};
    \end{array}
\right.
\end{eqnarray*}
\item[{\em(v)}] the sum $\sum_{k=1}^{\infty}a_k$ converges to $a < 1$.
\end{description}
\end{Lemma}
\Proof
If $p, q \in \natgtI$ are coprime then {\em(i)} certainly holds,
in which case {\em(ii)},
{\em(iii)} and {\em(iv)} uniquely determine
the $m_i$, $n_i$ and $a_k$. I claim that
for all sufficiently large rational $p$ and $q$,
the sum in {\em(v)} converges to a rational limit $a < 1$.
We may then take $p, q \in \natgtI$ large and coprime to complete the proof.

To prove the claim, apply standard facts about series of non-negative
terms to show that the sum $a = \sum_{k=1}^{\infty}a_k$ converges
for $p, q > 1$ and may be rearranged as follows:
$$
\begin{array}{rcc}
a &=&
  \sum_{i=0}^{\infty} (p^{-m_i}q^{-n_i} + 
    m_i p^{-m_i}q^{-n_i} +
    m_i n_i p^{-m_i}q^{-n_i} + 
    n_i p^{-m_i}q^{-n_i}) \\ \\
&=&
\left(
\begin{array}{r@{}c@{}l}
 (\sum_{i=0}^{\infty} p^{-m_i}q^{-n_i}) &{}+{}& 
    (\sum_{i=0}^{\infty} m_i p^{-m_i}q^{-n_i}) \\ &{}+{}& \\
 (\sum_{i=0}^{\infty} m_i n_i p^{-m_i}q^{-n_i}) &{}+{}& 
    (\sum_{i=0}^{\infty} n_i p^{-m_i}q^{-n_i})
\end{array}\right) \\ \\
&=&
\left(
\begin{array}{r@{}c@{}l}
(\sum_{m=2}^{\infty} p^{-m})(\sum_{n=2}^{\infty} q^{-n}) &{}+{}&
    (\sum_{m=2}^{\infty} mp^{-m})(\sum_{n=2}^{\infty} q^{-n} ) \\ &{}+{}& \\
(\sum_{m=2}^{\infty} mp^{-m})(\sum_{n=2}^{\infty} nq^{-n}) &{}+{}&
    (\sum_{m=2}^{\infty} p^{-m})(\sum_{n=2}^{\infty} nq^{-n}) \\
\end{array}\right) \\ \\
&=&
f(p)f(q) + g(p)f(q) + g(p)g(q) + f(p)g(q)
\end{array}
$$
where, using the formulas
$\sum_{i=0}^{\infty}x^i = (1-x)^{-1}$ and
$\sum_{i=1}^{\infty}ix^{i-1} = (1-x)^{-2}$
for the sums of a geometric series and its derivative, we have
$f(r) \IsDef (1 - r^{-1})^{-1} - 1 - r^{-1}$ and
$g(r) \IsDef r^{-1}[(1 - r^{-1})^{-2} - 1]$.
Thus $\sum_{k=1}^{\infty}a_k$ exists and
is a rational function, $a(p, q)$ say, with rational coefficients,
of the numbers $p$ and $q$ that tends to 0 as $p + q$ tends to $\infty$.
Hence for all sufficiently large rational $p$ and $q$, the
sum $\sum_{k=1}^{\infty}a_k$ converges to a rational $a = a(p, q)<1$.
(In fact, for $p, q \in \natgtI$, $a(p, q) < 1$
unless $p = q = 2$ or $\{p, q\} = \{2, 3\}$.)
\Done

We plan to encode the sequence $a_1, a_2, \ldots$ as the lengths of
line segments $[\Vv_1, \Vv_2]$, $[\Vv_2, \Vv_3]$, \ldots inscribed in the
unit circle of a 2-dimensional normed space $\sJ$.  So that the construction
works over an arbitrary field of scalars, we will arrange for the $\Vv_k$ to
have rational coordinates with respect to a basis $\Ve_1, \Ve_2$.  We
will also arrange for the norm of the vectors $\Vv_{k+1} - \Vv_{k}$ on $\sJ$
to agree with the 1-norm with respect to this basis. The following lemma will
give us gradients for the vectors $\Vv_{k+1} - \Vv_{k}$ that let
the $\Vv_k$ fit conveniently in the unit circle of $\sJ$.

\begin{Lemma}\label{lma:gradients}
Let $a_1, a_2, \ldots \in \ratgtO$ be as in lemma~\ref{lma:a-etc}.
Then there are $b_1, b_2, \ldots \in \ratgtO$
such that $1 > b_1 > b_2 > \ldots $ and
$b = \sum_{k=1}^{\infty} b_ka_k \in \ratgtO$.
\end{Lemma}
\Proof
Let $p$, $q$, etc.~be as in lemma~\ref{lma:a-etc}.  As a first approximation to
the $b_k$, let $b'_{4i+j} = p^{-m_i}q^{-n_i}$ for $i = 0, 1, \ldots$, $j = 1,
2, 3, 4$.  Then $b = \sum_{k=1}^{\infty}b'_ka_k = a(p^2, q^2) \in \ratgtO$
(where the rational function $a(r, s)$ is as in the proof of
lemma~\ref{lma:a-etc}). Now the sequence $b_1', b_2', \ldots$ is not strictly
decreasing, but we can derive a suitable strictly decreasing sequence $b_1,
b_2, \ldots $ from it. We do this in stages: at the $i$-th stage we construct a
block of 4 new values $b_{4i+1} > b_{4i+2} > b_{4i+3} > b_{4i+4}$ such that
$\sum_{j=1}^4b_{4i+j}a_{4i+j} = \sum_{j=1}^4b'_{4i+j}a_{4i+j}$ given a strict
upper bound $U_i$ for $b_{4i+1}$ and a strict lower bound $L_i$ for $b_{4i+4}$.
In stage 0, $U_0 = 1$.  Thereafter $U_i$ is the value $b_{4(i-1)+4}$ that ends
the block constructed in the previous stage.  We will arrange for $U_i >
p^{-m_i}q^{-n_i}$ (this is true for $i=0$ since $p^{-m_0}q^{-n_0} = (pq)^{-1} <
1$).  $L_i$ is $p^{-m_{i+1}}q^{-n_{i+1}}$ at the $i$-th stage for every $i$.
In the $i$-th stage, for $j = 1, 2, 3$, define $c_j(\delta) = b'_{4i+j} + (4 -
j)\delta = p^{-m_i}q^{-n_i} + (4 - j)\delta$ and define $c_4(\delta)$ to be the
rational function of $\delta$ that makes the following hold:
$$
\sum_{j=1}^4 c_j(\delta)a_{4i+j} = \sum_{j=1}^4 b'_{4i+j}a_{4i+j}.
$$
Then for $\delta \in
\ratgtO$, each $c_j(\delta) \in \rat$. Also $c_1(\delta) > c_2(\delta) >
c_3(\delta) > p^{-m_i}q^{-n_i} > c_4(\delta)$ and each $c_j(\delta)$ tends to
$p^{-m_i}q^{-n_i}$ as $\delta$ tends to 0.  So we may choose $\delta \in
\ratgtO$ such that with $b_{4i+j} = c_j(\delta)$ for $j = 1, 2, 3, 4$, we have
$U_i > b_{4i+1} > b_{4i+2} > b_{4i+3} > p^{-m_i}q^{-n_i} > b_{4i+4} > L_i$.  We
then have $U_{i+1} = b_{4i+4}
> L_i = p^{-m_{i+1}}q^{-n_{i+1}}$ so
the precondition for the next stage is satisfied.  Clearly, we have
$\sum_{k=1}^{4i+4}b_ka_k = \sum_{k=1}^{4i+4}b'_ka_k$ for all $i$, so that
$\sum_{k=1}^{\infty}b_ka_k = \sum_{k=1}^{\infty}b'_ka_k = b
\in \ratgtO$ and the sequence $b_k$ is as required.
\Done

We will now construct a 2-dimensional normed space $\sJ$ over an arbitrary ordered
field $K$ in which the graph of natural number multiplication is
definable.

\begin{figure}[h]
\begin{center}
\includegraphics[angle=0,scale=0.3]{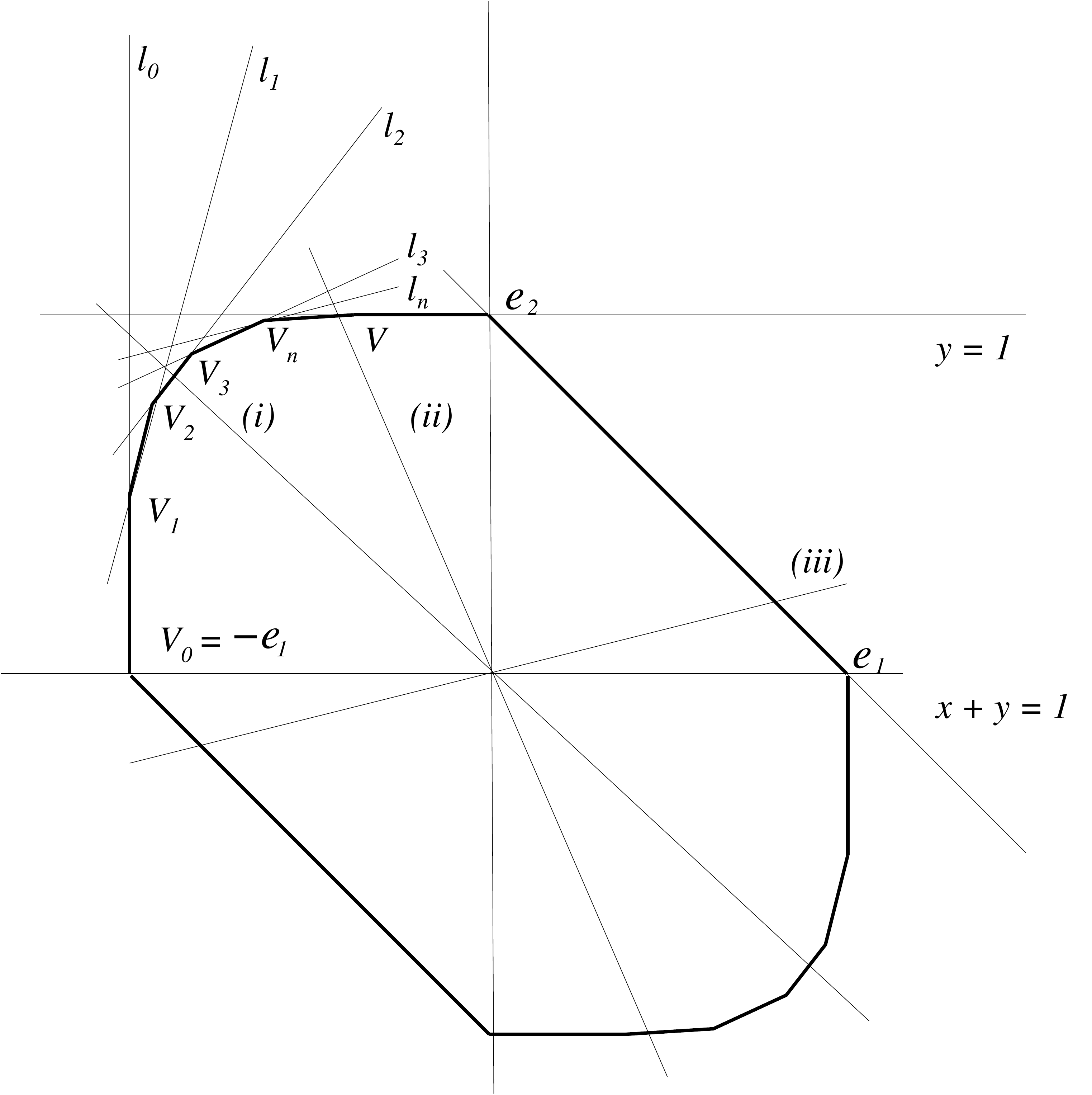}
\caption{The unit circle in the normed space $\sJ$}
\label{fig:J}
\end{center}
\end{figure}

Before embarking on the construction, note that the the basic ideas of affine
geometry and convexity theory all carry over to a vector space over an
arbitrary ordered field $K$. If $V$ is a vector space over $K$, we may define a
norm on $V$ by specifying its unit disc, which can be any convex subset $D$ of
$V$ that meets every line through the origin in a closed line segment $[-\Vx,
\Vx]$ where $\Vx \not= \VO$. Given such a $D$ and $\Vv \in V$, there is a
unique non-negative $r \in K$ such that $\Vv \in rD$ and $\Vv \not\in sD$ for $0 \le s < r$.
We set $\Norm{\Vv} = r$ and verify that this satisifies the norm axioms in
exactly the same way as when $K = \real$.  In examples over $\real$, the
requirements on $D$ are often simple consequences of the Heine-Borel theorem.
But that theorem  holds in $K$ iff $K$ is isomorphic to $\real$, so this method
of proof does not apply in general. With $K = \rat$, for example, if we
take $D = \{(x, y) \in \rat^2 \ST x^2 + y^2 \le 1\}$, then $D$ is
closed, convex and bounded, but with $\Vv = (1, 1)$, the
set of $s$ such that $\Vv \not\in sD$ is non-empty and bounded but has no least upper bound.

%However, a bounded and
%topologically closed 1-dimensional convex set need not be a closed line segment $[\Vp, \Vq]$.
%Indeed, any ordered field that is not isomorphic to $\real$ is
%order-incomplete: i.e., there is a bounded sequence
%$s_1 < s_2 < \ldots < s$ say, that has no least upper bound.
%Then $\{x \ST \ex{i} s_1 \le x < s_i\}$ is convex bounded, topologically
%closed but not equal to $[s_1, t]$ for any $t$.
%Thus we will often be obliged to prove the condition on the intersection
%on lines with $D$ algebraically rather than by a compactness argument.

Let $\Ve_1$ and $\Ve_2$ be the standard basis for $K^2$.
Define vectors $\Vv_k$, $k = 0, 1 \ldots$
in the north-west quadrant as follows:
\begin{eqnarray*}
\Vv_0 &\IsDef& -\Ve_1 \\
\Vv_1 &\IsDef& -\Ve_1 + (1 - b)\Ve_2 \\
\Vv_{k+1} &\IsDef& \Vv_k + (1 - b_k)a_k\Ve_1 + b_ka_k\Ve_2
\end{eqnarray*}
where the $a_k$, $b_k$ and $b$ are the rational numbers of lemmas~\ref{lma:a-etc} and~\ref{lma:gradients}.
See Figure~\ref{fig:J}.

If $l$ is any line in $K^2$, we define its {\em
gradient} to be the symbol $\infty$ if $l$ is parallel to $\Ve_2$
or to be the unique $g \in K$ such that $l$ is parallel to $\Ve_1 + g\Ve_2$
otherwise.

Let the line $l_k$ through $\Vv_{k}$ and $\Vv_{k+1}$ have gradient $g_k$.
So $g_0 = \infty$, and $g_k = b_k/(1-b_k) = (b_k^{-1}-1)^{-1}$, for $k = 1, 2,
\ldots$  By lemma~\ref{lma:gradients}, $1 > b_1 > b_2 > \ldots > 0$, so that $g_1 >
g_2 > \ldots > 0$.  Also, if we write $\Vv_k = x_k\Ve_1 + y_k\Ve_2$, the
sequence $(x_k, y_k)$ tends to a limit $(a - b - 1, 1)$ in $\rat^2$ as $k$
tends to $\infty$, where $a$ is as in lemma~\ref{lma:a-etc} so that $-1 < a -
b - 1 < 0$.

Define a subset $D$ of $K^2$ as follows:
$$
D \IsDef
   (\bigcap_{k=0}^{\infty}H_k) \cap A \cap B \cap
   (\bigcap_{k=0}^{\infty}-H_k) \cap -A \cap -B
$$
where $H_k$ is the closed half-plane that contains the origin and has
the line $l_k$ as boundary and where $A$ and $B$ are the closed half-planes 
defined by the formulas $y \le 1$ and $x + y \le 1$ respectively.
Clearly $D$ is convex and symmetric about the origin.

I claim that $D$ meets every
line through the origin in a line segment $[-\Vx, \Vx]$ where $\Vx \not= \VO$.
To prove this, first note that by routine algebra,
if $h_k = y_k/x_k$ is the gradient of the line through $\VO$ and $\Vv_k$, then
$0 = h_0 > h_1 > h_2 > \ldots$ Also from the remarks above about the $x_k$ and
$y_k$, the $h_k$ tend in $\rat$ to the limit $h = (a - b - 1)^{-1}$.
Now let $l$ be a line through the origin with gradient $g$. We need
to exhibit an $\Vx \not= \VO$ such that $l \cap D = [-\Vx, \Vx]$.
By symmetry, it is enough to find $\Vx \not= \VO$ such that $l \cap D_u = [\VO, \Vx]$
where $D_u = D \cap H_u$, $H_u$ being the half-plane defined by the formula
$y \ge 0$.
We identify three cases as follows (see Figure~\ref{fig:J}):

{\em(i) $h_k > g \ge h_{k+1}$ for some $k$:}
clearly the $\Vv_i$ lie in the interior of the half-planes $A$ and $B$;
also, the conditions on the gradients $g_j$ imply that $\Vv_i \in H_j$
for all $i, j$, and so as $D$ is convex, $[\Vv_k, \Vv_{k+1}] \subseteq
D$.  As $h_k > g \ge h_{k+1}$, $l$ meets $[\Vv_k, \Vv_{k+1}]$ at
some point $\Vx \not= \VO$ and then as every neighbourhood of $\Vx$ meets
both $D$ and its complement, we must have $l \cap D_u = [\VO, \Vx]$.

{\em(ii) $g \not= \infty$ and $h_k > g$ for all $k$:}
in this case, $l$ meets the line $y = 1$
at a point $\Vx = p\Ve_1 + \Ve_2 \not= \VO$ where either $a - b - 1 \le p < 0$
or $(a - b - 1) - p$ is a positive infinitesimal.
$\Vx$ then lies on the boundary of half-plane $A$ and in the interior of half-plane $B$.
Since the gradients $g_k = (y_{k+1}-y_k)/(x_{k+1}-x_k)$ are rational and
satisfy $g_1 > g_2 > \ldots > 0$ and since the
points $(x_k, y_k)$ converge in $\rat^2$ to the limit $(a-b-1, 1)$,
the line $l_k$ must meet the line $y = 1$ at a point
$r\Ve_1 + \Ve_2$ where $r < a - b - 1$ is rational. But this implies
that $r < p$ so that $l_k$ and $l$ meet at a point
$\Vy = s\Ve_1 + t\Ve_2$ where $r < s < p$ and $t > 1$
so that $\Vx$ is to the south and east of the point $\Vy \in l_k$.
Thus $\Vx \in H_k$ for every $k$ and every neighbourhood
of $\Vx$ meets both $D$ and its complement, whence $l \cap D_u = [\VO, \Vx]$. 

{\em(iii) $g \ge 0$ or $g = \infty$:}
it is easy to see that the intersection of $D$ with the
north-east quadrant is the triangle $\triangle \VO\Ve_1\Ve_2$, so that,
if $g \ge 0$ or $g = \infty$, $l$ meets $D_u$ in the
interval $[\VO, \Vx]$, where $\Vx = \Ve_2$ if $g = \infty$ and
$\Vx = (\Ve_1 + g\Ve_2)/(1 + g)$ otherwise. In both cases $\Vx \not= \VO$.

We have proved that $D$ is convex and meets each 
line through the origin in a line segment $[-\Vx, \Vx]$ where $\Vx \not= \VO$.
Hence there is a norm on $K^2$ having $D$ as its unit disc.
Define $\sJ$ to be $K^2$ equipped with that norm.

The case analysis on the gradient $g$ of the line $l$ in the argument
above shows that the upper half of the unit circle in
$\sJ$ comprises:
{\em(i)} the line segments $[\Vv_0, \Vv_1]$, $[\Vv_1, \Vv_2]$ \ldots,
{\em(ii)} the set $E$ of all points on the line segment
$[-\Ve_1 + \Ve_2, \Ve_2]$ that lie to the east of every $\Vv_k$
and {\em(iii)} the line segment $[\Ve_1, \Ve_2]$.
If $K$ is archimedean, $E$ is the line segment $[\Vv, \Ve_2]$
where $\Vv = (a-b-1)\Ve_1 + \Ve_2$, but in the non-archimedean case, $E$
also contains every point $\Vv - \epsilon_*\Ve_1$ where $\epsilon_*$
is a positive infinitesimal.

Having defined $\sJ$ and described its unit circle, let us develop
the formula $\pM{x, y, z}$.
First we define:
\begin{eqnarray*}
\pEP{\Vp} &\IsDef& \all{\Vu\;\Vw} \Norm{\Vu} = \Norm{\Vp} = \Norm{\Vw} \And \Vp = \frac{1}{2}(\Vu + \Vw) \Imp \Vu = \Vp = \Vw\\
\pSEP{\Vp} &\IsDef& \pEP{\Vp} \And \ex{\Vu\;\Vw} \pEP{\Vu} \And \pEP{\Vw} \And \Norm{\Vu} = \Norm{\Vw} = \Norm{\Vp} \And {} \\
 & & \quad\quad\quad\quad\quad\quad\quad
     0 \not= \Norm{\Vp - \Vu} \not= \Norm{\Vp - \Vw} \not = 0 \\
\pADS{\Vp, \Vq} &\IsDef& \pSEP{\Vp} \And \pSEP{\Vq} \And \Vp \not= \Vq \And \Norm{\Vp} = \Norm{\Vq} = \Norm{(\Vp + \Vq)/2}.
\end{eqnarray*}
Thus in any normed space $\pEP{\Vp}$ holds iff $\Vp$ is an extreme
point of the sphere $S_{\Norm{\Vp}}$ of radius $\Norm{\Vp}$ centred on the origin and $\pSEP{\Vp}$ holds iff $\Vp$ is an extreme point of $S_{\Norm{\Vp}}$ that
is not equidistant from every other extreme point of $S_{\Norm{\Vp}}$.
Let us call points $\Vp$ satisfying $\pSEP{\Vp}$ {\em special} extreme points
of $S_{\Norm{\Vp}}$.
In $\sJ$, the special extreme points are just the non-zero ones.
In any normed space,
$\pADS{\Vp, \Vq}$ holds iff $\Vp, \Vq$ are adjacent special extreme points of $S_{\Norm{\Vp}}$.
Next we define:
\begin{eqnarray*}
\pHPV{\Vp_1, \ldots, \Vp_n} &\IsDef&
   \Ands_{i=1}^{n-1} \pADS{\Vp_i, \Vp_{i+1}} \And
   \Ands_{i=1}^{n-1}\Ands_{j=i+1}^{n} \Vp_i \not= \Vp_j \\
\pHPL{x_1, \ldots, x_n} &\IsDef&
   \ex{\Vp_1\;\ldots\Vp_{n+1}} \pHPV{\Vp_1, \ldots, \Vp_{n+1}} \And {}\\
 & &	\Ands_{i=1}^n x_i = \Norm{\Vp_{i+1} - \Vp_{i}}.
\end{eqnarray*}
So in any normed space, $\pHPV{\Vp_1, \ldots, \Vp_n}$ holds iff $\Vp_1, \ldots, \Vp_n$ is
a sequence of special extreme points of $S_{\Norm{\Vp_1}}$
forming the vertices of a Hamiltonian path made up of straight line
segments inscribed in $S_{\Norm{\Vp_1}}$.
$\pHPL{x_1, \ldots, x_n}$ holds
iff $x_1, \ldots, x_n$ are the lengths of the successive edges of such a path.
Finally we define:
\begin{eqnarray*}
\pMGI{x, y,z} &\IsDef& \ex{u} \pHPL{1, x, z, y, u} \And
       1 < x < z > y > u < 1 \\
\pM{x, y, z} &\IsDef& \pMGI{x+2, y+2, 4 + 2x + 2y + z}.
\end{eqnarray*}

Identifying $\nat$ and $\nat_K$, I claim that $\pMGI{x, y, z}$ holds in $\sJ$
iff $x, y, z \in \natgtI$ and $z = xy$. This implies that $\pM{x, y , z}$
defines the graph of multiplication in $\nat$.

To prove the claim,
first note that for vectors in the north-east quadrant, the $\sJ$-norm
is the same as the 1-norm, so $\Norm{\Vv_{k+1} - \Vv_k}_{\sJ} = (1 - b_k)a_k +
b_ka_k = a_k$, for $k > 0$. Also $\Norm{\Ve_2 - \Ve_1}_{\sJ} > 1$ because
$\Ve_2 - \Ve_1 \not\in D$.
So $\Norm{\Ve_2 - \Ve_1}_{\sJ}$ is greater than the $\sJ$-length of any line segment
that can be inscribed in the north-west quadrant of the unit circle of $\sJ$.

Now, if $\pMGI{x, y, z}$ holds, there are special extreme points $\Vp_1,
\ldots, \Vp_6$ forming the vertices of a Hamiltonian path inscribed in $S_r$
for some $r > 0$ such that the edges $[\Vp_2, \Vp_1]$, \ldots, $[\Vp_6, \Vp_5]$
have $\sJ$-lengths $1, x, y, z, u$ respectively where $1 < x < z > y > u < 1$ for
some $u$.  But then the local maximum $z$ can only be the $\sJ$-length of
${\pm}[r\Vv_{4i+3}, r\Vv_{4i+4}]$ for some $i$ (it cannot be $\Norm{r\Ve_2 -
r\Ve_1}$ since, even when $K$ is archimedean, there is no Hamiltonian path in $S_r$ with 6 extreme
points as vertices such that
${\pm}[r\Ve_1, r\Ve_2]$ is the $3^{\mbox{\scriptsize rd}}$ edge). As $u < 1$,
$\Vp_1, \ldots, \Vp_6$, are ${\pm}\Vv_{4i+1}, \ldots, {\pm}\Vv_{4i+6}$ in that
order, so that
$1 = \Norm{\Vp_2 - \Vp_1} = ra_{4i+1} = rp^{-m_i}q^{-n_i}$ and $r =
p^{m_i}q^{n_i}$. Hence $x = rm_ip^{-m_i}q^{-n_i} = m_i$ and, similarly,
$z = m_in_i$ and $y = n_i$, so $x, y, z \in \natgtI$ and $z = xy$.

Conversely, if $x, y, z \in \natgtI$ and $z = xy$, then for some $i$,
$x = m_i$ and $y = n_i$. Let $r = p^{m_i}q^{n_i}$ so that $r\Vv_{4i+1}, \ldots,
r\Vv_{4i+6}$ are the vertices of a Hamiltonian path inscribed in $S_r$, whose
edges have $\sJ$-lengths $1, x, z, y, rp^{-m_{i+1}}q^{-n_{i+1}}$
respectively. 
By lemma~\ref{lma:a-etc},
we have $1 < x < z > y > rp^{-m_{i+1}}q^{-n_{i+1}} < 1$,
and so $\pMGI{x, y, z}$
holds in $\sJ$ with $rp^{-m_{i+1}}q^{-n_{i+1}}$ as the witness for $u$.

\begin{Theorem}\label{thm:main-result}
There is formula $\pM{x, y, z}$ in the purely additive
language of normed spaces $\LNA$ such that for any ordered field $K$ and any $d
\in \{2, 3, 4, \ldots\} \cup \{\infty\}$, there is a normed
space $\sJ^d$ of dimension $d$ in which $\pM{x, y, z}$ defines
the graph of natural number multplication on $\nat_K$.  
\end{Theorem}
\Proof
Let $\pM{x, y, z}$ and $\sJ$ be as above. Let $W$ be a $(d-2)$-dimensional
vector space over $K$ equipped with the 1-norm with respect to some basis $\Vb_1, \Vb_2, \ldots$
($\Norm{\sum_i c_i\Vb_{i}}_W = \sum_i|c_i|$)
and let $\sJ^d$ be the 1-sum $\sJ \times W$ ($\Norm{(\Vv, \Vw)}_{\sJ^d} =
\Norm{\Vv}_{\sJ} + \Norm{\Vw}_W$). Identify $\sJ$ and $W$ with $\sJ \times 0$ and $0
\times W$ respectively.
Then the extreme points of $S_r$ in
$\sJ^d$ comprise the extreme points of $S_r \cap \sJ$ together with the
extreme points ${\pm} r\Vb_1, {\pm} r\Vb_2, \ldots$ of $S_r \cap W$.
Moreover, $\Norm{{\pm} r\Vb_n - \Vp} = 2r$ for every
extreme point $\Vp$ of $S_r$ with $\Vp \not= {\pm} r\Vb_n$. Thus, in the sense
defined above, the only special extreme points of $S_r$ in $\sJ^d$ are those of
$\sJ$. It follows from the discussion above that $\pM{x, y, z}$ defines the
graph of natural number multiplication in $\sJ^d$.
\Done

\begin{Corollary}\label{cor:arb-ns-undec1}
If $\cal C$ is any non-empty class of ordered fields, then the theories
$\NSA{\relax}({\cal C})$,
$\NSA{n}({\cal C})$ for $1 < n \in \nat$, 
$\NSA{\fin}({\cal C})$
and $\NSA{\infty}({\cal C})$
are all undecidable.
\end{Corollary}
\Proof
This is immediate from theorems~\ref{thm:gen-linear-undec} and~\ref{thm:main-result}.
\Done

\begin{Corollary}\label{cor:arb-ns-undec2}
If $\cal C$ is any non-empty class of ordered fields, then the theories
$\NS{\relax}({\cal C})$,
$\NS{n}({\cal C})$ for $1 < n \in \nat$, 
$\NS{\fin}({\cal C})$
and $\NS{\infty}({\cal C})$
are all undecidable.
\end{Corollary}
\Proof
By the corollary, we cannot decide sentences in the additive fragments
of these theories.
\Done

When $\cal C$ is any class of ordered fields including
$\rat$, the theory of $\cal C$ is then undecidable by a classic result of Julia
Robinson \cite{robinson49}. Corollary~\ref{cor:arb-ns-undec2} tells us nothing
new for such a $\cal C$, but as the additive fragment of the theory
of ordered fields is decidable, Corollary~\ref{cor:arb-ns-undec1} is significant.
Both corollaries have force when the theory of the fields in $\cal C$ is
decidable, e.g., when the fields in $\cal C$ are real closed.

%For completeness, given the title of the present paper, I should mention the
%case when $\cal C$ is empty. The proof that the various theories are all
%decidable in that case can safely be left to the reader.

% \newpage

\bibliographystyle{plain}

\bibliography{bookspapers}

\end{document}